\input amstex 
\documentstyle{amsppt}
\document
\define\m{\medskip}
\define\C#1{\Cal #1}
\define\p{{\m \noindent \it Proof. }}
\def\Im{\operatorname{Im}\,}

\def\hr{\operatorname{hr}}
\topmatter
\title
On certain geometric and  homotopy properties of closed symplectic manifolds
\endtitle
\author
Ra\'ul Ib\'a\~nez, Yuli Rudyak, Aleksy Tralle and Luis Ugarte
\endauthor
\address
Departamento de Matem\'aticas, Facultad de
Ciencias, Universidad del Pais Vasco, Apdo. 644, 48080 Bilbao, Spain,
\endaddress
\vskip6pt
\email:  mtpibtor\@lg.ehu.es
\endemail
\address
FB6/Mathematik, Universit\"at Siegen, 57068 Siegen, Germany
\endaddress
\vskip6pt
\email
rudyak\@mathematik.uni-siegen.de,
july\@mathi.uni-heidelberg.de
\endemail
\vskip6pt
\address
 University of Warmia and Masuria, 10561 Olsztyn, Poland
\endaddress
\vskip6pt
\email
tralle\@tufi.wsp.olsztyn.pl
\endemail
\vskip6pt
\address Departamento de Matem\'aticas, Facultad de
Ciencias, Universidad de Zaragoza, 50009 Zaragoza, Spain
\endaddress
\vskip6pt
\email: ugarte\@posta.unizar.es
\endemail
\endtopmatter
\rightheadtext{Homotopy properties of symplectic manifolds}

\head 1. Introduction
\endhead

Homotopy properties of closed symplectic manifolds attract the
attention of geometers since the classical papers of Sullivan \cite{S}
and Thurston \cite{Th}. On one hand, "soft" homotopy techniques help
in the solution of many "hard" problems in symplectic geometry, cf.
\cite{G1, McD, RT1, TO}. On the other hand, it is still unknown if
there are specific homotopy properties of closed  manifolds dependent
on the existence of  symplectic structures on them. It turns out that
symplectic manifolds violate many specific homotopy conditions shared
by the K\"ahler manifolds (which form a subclass of symplectic
manifolds). In particular, if $M$ is a closed K\"ahler manifold then
the following holds:
\roster
\item all the odd-degree Betti numbers $b_{2i+1}(M)$ are even;
\item $M$ has the Hard Lefschetz property;
\item all Massey products (of all orders) in $M$ vanish.
\endroster
It is well known (and we shall see it below) that closed symplectic
manifolds violate {\it all} the homotopy properties (1) -- (3).
However, it is not clear whether properties (1) -- (3) are {\it
independent} or not, in case of closed symplectic manifolds or certain
classes of these ones. In other words, can a combination of the type

\centerline{(1) -- (2) -- non-(3)}
\noindent
 be realized by a closed symplectic manifold (possibly, with prescribed properties).
 The knowledge of an answer to this question might shed a new light on
the whole understanding of closed symplectic manifolds.
\m
In Theorem 3.1 we have summarized our knowledge by writing down the
corresponding tables. We have considered two classes of symplectic
manifolds: the class of symplectically aspherical symplectic manifolds
and the class of simply-connected symplectic manifolds. Recall that a
symplectically aspherical manifold is a symplectic manifold
$(M,\omega)$ such that $\omega|\pi_2(M)=0$, i.e.
$$
\int_{S^2}f^{\#}\omega=0
$$
for every map $f: S^2\to M$. In view of the Hurewicz Theorem, a closed
symplectically aspherical manifold always has a non-trivial
fundamental group. It is well known that symplectically aspherical
manifolds play an important role in geometry and topology of
symplectic manifolds, \cite{F, G2, H, RO, RT2}.

The next topic of the paper is about symplectically harmonic forms on
closed symplectic manifolds. Brylinski \cite{B} and Libermann (Thesis,
see~\cite{LM}) have introduced the concept of a {\it symplectic star
operator} $*$ on a symplectic manifold. In a sense, it is a symplectic
analog of the Hodge star operator which is defined in terms of the
given symplectic structure $\omega$. Using this operator, one defines
a symplectic codifferential $\delta:=(-1)^{k+1}(* d*),\, \deg\,\delta
=-1$. Now we define {\it symplectically harmonic} differential forms
$\alpha$ by the condition
$$
\delta\alpha=0,\,\,\,\, d\alpha=0.
$$
Let $\Omega_{\text{hr}}^*(M,\omega)$ denote the space of all
symplectically harmonic forms on $M$. Clearly, the space
$H^k_{\hr}(M):=\Omega^k_{\hr}/(\Omega^k_{\hr}\cap \Im d)$  is a
subspace of the de Rham cohomology space $H^k(M)$.
\m
Here we also have an interesting relation between geometry and
homotopy theory. For example, Mathieu~\cite{M} proved that $
H^k_{\hr}(M, \omega)= H^k(M)$ if and only if $M$ has the Hard
Lefschetz property. We will also see that the Lefschetz map
$$
L^k:H^{m-k}(M) \to H^{m+k}(M),\ \dim\,M=2m
$$
(multiplication by $[\omega]^k$) plays an important role in studying
of $H^k_{\hr}(M, \omega)$.

\m We set $h_k(M,\omega)=\dim\, H^k_{\hr}(M, \omega)$. According to
Yan~\cite{Y}, the following question was posed by Boris Khesin and
Dusa McDuff.
\m
{\bf Question:} Are there closed manifolds endowed with a continuous
family $\omega_t$ of symplectic structures such that $h_k(M,
\omega_t)$ varies with respect to $t$?
\m
Yan \cite{Y} constructed a closed 4-dimensional manifold $M$ with
varying $h_3(M)$. So, he answered affirmatively the above question.

(Actually,  Proposition 4.1 from \cite{Y} is wrong, the
Kodaira--Thurston manifold is a counterexample, but its Corollary 4.2
from \cite{Y} is correct because it follows from our Lemma 4.4. Hence,
the whole construction holds.)
\m
However, the Yan's proof was essentially 4-dimensional. Indeed,
Yan~\cite{Y} wrote:

``For higher dimensional closed symplectic manifolds, it is not clear
how to answer the question in the beginning of this section",

i.e. the above stated question.
\m
In this note we prove the following result (Theorem 4.6): There exists
at least one 6-dimensional indecomposable closed symplectic manifold
$N$ with varying $h_5(N)$.
\m
Moreover, Yan remarked that there is no 4-dimensional closed
symplectic {\it nilmanifolds} $M$ with varying $\dim\,H^*_{\hr}(M)$.
On the contrary, our example is a certain 6-dimensional nilmanifold.

\head 2. Preliminaries and notation
\endhead

Given a topological space $X$, let $(\C M_X,d)$ be the Sullivan
model of $X$, that is, a certain natural commutative DGA algebra over
the field of rational numbers $\Bbb Q$ which is a homotopy invariant
of $X$, see \cite{DGMS, TO, S} for details. Furthermore, if $X$ is a
nilpotent $CW$-space of finite type then $(\C M_X,d)$ completely
determines the rational homotopy type of $X$.
\m
A space $X$ is called {\it formal} if there exists a DGA-morphism
$$
\rho: (\C M_X,d)\to (H^*(X;\Bbb Q),0)
$$
inducing isomorphism on the cohomology level. Recall that every closed
K\"ahler manifold in formal~\cite{DGMS}.
\m We refer the reader to \cite{K, Ma, RT1} for the definition
of Massey products. It is well known and easy to see that Massey
products yield an obstruction to formality \cite{DGMS, RT1, TO}. In
other words, if the space is formal then all Massey products must be
trivial. Thus, all the Massey products in every K\"ahler manifold
vanish.
\m We need also the following result of Miller~\cite{Mi}:
\proclaim{2.1. Theorem} Every closed simply-connected manifold $M$ of
dimension $\leq 6$ is formal. In particular, all Massey products in $M$ vanish.
\qed
\endproclaim
\m
The next homotopy property related to symplectic (in particular,
K\"ahler) structures is the {\it Hard Lefschetz property}. Given a
symplectic manifold $(M^{2m},\omega)$, we denote by $[\omega]\in
H^2(M)$ the de Rham cohomology class of $\omega$. Furthermore, we
denote by $L_{\omega}: \Omega^k(M)\to\Omega^{k+2}(M)$ the
multiplication by $\omega$ and by $L_{[\omega]}: H^k(M)\to H^{k+2}(M)$
the induced homomorphism in the de Rham cohomology $H^*(M)$. As usual
we write $L$ instead of $L_{\omega}$ or $L_{[\omega]}$ if there is no
danger of confusion. We say that a symplectic manifold
$(M^{2m},\omega)$ has the {\it Hard Lefschetz property} if, for every
$k$, the homomorphism
$$
L^k: H^{m-k}(M)\to H^{m+k}(M)
$$
is surjective. In view of the Poincar\'e duality, for closed manifolds
$M$ it means that every $L^k$ is an isomorphism. We need also the
following result of Gompf~\cite{G1, Theorem 7.1}.

\proclaim{2.2. Theorem} For any even dimension $n\geq 6$, finitely
presented group $G$ and integer $b$ there is a closed symplectic
$n$-manifold $M$ with $\pi_1(M)\cong G$ and $b_i(M)\geq b$ for $2\leq
i\leq n-2$, such that $M$ does not satisfy the Hard Lefschetz
condition. Furthermore, if $b_1(G)$ is even then all degree-odd Betti
numbers of $M$ are even.
\qed
\endproclaim

We denote such manifold $M$ by $M(n,G,b)$.
\m
{\noindent \bf 2.3. Remark.} Theorem 7.1 in \cite{G1} is formulated in
a slightly different way, but the {\it proof} is based on constructing
of $M$ by some "symplectic summation" in a way to violate the Hard
Lefschetz property.
\m
In our explicit constructions  we will need some particular classes of
manifolds, namely, {\it nilmanifolds}, resp. {\it solvmanifolds}.
These are homogeneous spaces of the form $G/\Gamma$, were $G$ is a
simply connected nilpotent, resp. solvable Lie group and $\Gamma$ is a
co-compact discrete subgroup (i.e. a {\it lattice}). The most
important information for us is the following (see e.g. \cite{TO} for
the proofs):
\m
{\noindent \bf 2.4. Recollection.} {\rm (i)} Let $\frak g$ be a
nilpotent Lie algebra with structural constants $c_k^{ij}$ with
respect to some basis, and let $\{\alpha_1,...,\alpha_n\}$ be the dual
basis of $\frak g^*$. Then the differential in the
Chevalley--Eilenberg complex $(\Lambda^*\frak g^*,d)$ is given by the
formula
$$
d\alpha_k=-\sum_{1\leq i<j<k}c_k^{ij}\alpha_i\wedge\alpha_j.
$$
\par {\rm (ii)} Let $\frak g$ be the Lie algebra of a simply connected
nilpotent Lie group $G$. Then, by Malcev's theorem, $G$ admits a
lattice if and only if $\frak g$ admits a basis such that all the
structural constants are rational. Moreover, this lattice is unique up
to an automorphism of $G$.
\par {\rm (iii)} Let $G$ and $\frak g$ be as in (ii), and suppose
that $G$ admits a lattice $\Gamma$. By Nomizu's theorem, the
Chevalley--Eilenberg complex $(\Lambda^*\frak g^*,d)$ is
quasi-isomorphic to the de Rham complex of $G/\Gamma$. Moreover,
$(\Lambda^*\frak g^*,d)$ is a minimal differential algebra, and hence
it is isomorphic to the minimal model of $G/\Gamma$:
$$
(\Lambda^*\frak g^*,d)\cong (\C M_{G/\Gamma},d).
$$
Also, any cohomology class $[a]\in H^k(G/\Gamma)$ contains a
homogeneous representative $\alpha$. Here we call the form $\alpha$
{\it homogeneous} if the pullback of $\alpha$ to $G$ is left
invariant.
\m
Let $\omega_0$ be the standard symplectic form on $\Bbb CP^m$. Recall
that every closed symplectic manifold $(M^{2n},\omega)$ with integral
form $\omega$ can be symplectically embedded into $\Bbb CP^m$ for $m$
large enough, with the (known) smallest possible value of $m$ equal to
$n(n+1)$ \cite{Gr, Ti}. We will use the {\it blow-up} construction
with respect to such embedding \cite{McD, RT1}. We need the following
result.

\proclaim{2.5. Theorem} Let $(M^{2n},\omega)$ be a closed connected symplectic
manifold, let $i: (M,\omega) \to (\Bbb CP^m,\omega_0)$ be a symplectic
embedding, and let $\widetilde{\Bbb CP^m}$ be the blow-up along $i$. Then the following holds:
\par{\rm(i)}  $\widetilde{\Bbb CP^m}$ is a simply-connected symplectic manifold;
\par{\rm(ii)} if there exists $i$ such that $b_{2i+1}(M)$ is odd, then there
exists $k$ such that $b_{2k+1}(\widetilde{\Bbb CP^m})$ is odd;
\par{\rm(iii)} if $M$ possesses a non-trivial Massey triple product
and $m-n\geq 4$, then $\widetilde{\Bbb CP^m}$ possesses a non-trivial
Massey triple product. Moreover, if there is a non-trivial Massey
product $\langle\alpha, \beta, [\omega]\rangle\in H^*(M)$, $\alpha,
\beta\in H^*(M)$, then $\widetilde{\Bbb CP^m}$ possesses a non-trivial
Massey triple product even for $m-n=3$.
\endproclaim
\p (i) and (ii) are proved in \cite{McD}, (i) and (iii) are proved in
\cite{RT1}.
\qed

\head 3. Relation between homotopy properties of closed symplectic manifolds \endhead
\proclaim{3.1. Theorem}
The relations between the Hard Lefschetz property, evenness of
odd-degree Betti numbers and vanishing of the Massey products for
closed symplectic manifolds are given by the following tables:

                    {\rm TABLE 1:} symplectically aspherical case;

                    {\rm TABLE 2:} simply-connected case.
\endproclaim
\m
The word {\bf Impossible} in the table means that there is no closed
symplectic manifold (aspherical or simply connected) that realizes the
combination in the corresponding line.
\m
The sign {\bf ?} means that we (the authors) do not know whether a
manifold with corresponding properties exists.
\vfill
\pagebreak
\centerline{\bf Table 1: Symplectically Aspherical Symplectic Manifolds}
\vskip 0.5cm
\vbox{\tabskip=0pt\offinterlineskip
\def\tablerule{\noalign{\hrule}}
\halign to 350pt{\strut#&\vrule#\tabskip=1em plus2em&
\hfil# \hfil & \vrule#&
\hfil# \hfil & \vrule#&
\hfil# \hfil & \vrule#&
\hfil { #}  \hfil & \vrule # \tabskip=0pt\cr \tablerule
\omit&height5pt&\omit&height5pt&\omit&height5pt&\omit&height5pt&&\cr
&& Triviality of &&
\hidewidth Hard Lefschetz
\hidewidth &&
\hidewidth Evenness of \hidewidth && \hidewidth \phantom{K\"ahler ($\Bbb T^{2n}$)}
\hidewidth &\cr &&
\hidewidth Massey Products \hidewidth && Property &&
 $b_{2i+1}$  && \phantom{K\"ahler ($\Bbb T^{2n}$)}&\cr \tablerule
\omit&height5pt&\omit&height5pt&\omit&height5pt&\omit&height5pt&&\cr
&& yes && yes && yes && K\"ahler ($\Bbb T^{2n}$)&\cr
\omit&height5pt&\omit&height5pt&\omit&height5pt&\omit&height5pt&&\cr\tablerule
\omit&height5pt&\omit&height5pt&\omit&height5pt&\omit&height5pt&&\cr
&& yes && yes && no && Impossible &\cr
\omit&height5pt&\omit&height5pt&\omit&height5pt&\omit&height5pt&&\cr\tablerule
\omit&height5pt&\omit&height5pt&\omit&height5pt&\omit&height5pt&&\cr
&& yes && no && yes && ? &\cr
\omit&height5pt&\omit&height5pt&\omit&height5pt&\omit&height5pt&&\cr\tablerule
\omit&height5pt&\omit&height5pt&\omit&height5pt&\omit&height5pt&&\cr
&& yes && no && no && ? &\cr
\omit&height5pt&\omit&height5pt&\omit&height5pt&\omit&height5pt&&\cr\tablerule
\omit&height5pt&\omit&height5pt&\omit&height5pt&\omit&height5pt&&\cr
&& no && yes && yes && ? &\cr
\omit&height5pt&\omit&height5pt&\omit&height5pt&\omit&height5pt&&\cr\tablerule
\omit&height5pt&\omit&height5pt&\omit&height5pt&\omit&height5pt&&\cr
&& no && yes && no && Impossible &\cr
\omit&height5pt&\omit&height5pt&\omit&height5pt&\omit&height5pt&&\cr\tablerule
\omit&height5pt&\omit&height5pt&\omit&height5pt&\omit&height5pt&&\cr
&& no && no && yes && $K\times K$ &\cr
\omit&height5pt&\omit&height5pt&\omit&height5pt&\omit&height5pt&&\cr\tablerule
\omit&height5pt&\omit&height5pt&\omit&height5pt&\omit&height5pt&&\cr
&& no && no && no && $K$ &\cr
\omit&height5pt&\omit&height5pt&\omit&height5pt&\omit&height5pt&&\cr
\tablerule\noalign{\smallskip}\hfil\cr}}

\centerline{\bf Table 2: Simply-Connected Symplectic Manifolds}
\vskip 0.5cm
\vbox{\tabskip=0pt\offinterlineskip
\def\tablerule{\noalign{\hrule}}
\halign to 350pt{\strut#&\vrule#\tabskip=1em plus2em&
\hfil# \hfil & \vrule#&
\hfil# \hfil & \vrule#&
\hfil# \hfil & \vrule#&
\hfil# \hfil & \vrule#\tabskip=0pt\cr \tablerule
\omit&height5pt&\omit&height5pt&\omit&height5pt&\omit&height5pt&&\cr
&& Triviality of &&
\hidewidth Hard Lefschetz
\hidewidth &&
\hidewidth Evenness of \hidewidth && \hidewidth \phantom{Impossible} \hidewidth &\cr &&
\hidewidth Massey Products \hidewidth && Property &&
$b_{2i+1}$ && \phantom{Impossible}&\cr \tablerule
\omit&height5pt&\omit&height5pt&\omit&height5pt&\omit&height5pt&&\cr
&& yes && yes && yes && K\"ahler ($\Bbb C P^n$)&\cr
\omit&height5pt&\omit&height5pt&\omit&height5pt&\omit&height5pt&&\cr\tablerule
\omit&height5pt&\omit&height5pt&\omit&height5pt&\omit&height5pt&&\cr
&& yes && yes && no && Impossible &\cr
\omit&height5pt&\omit&height5pt&\omit&height5pt&\omit&height5pt&&\cr\tablerule
\omit&height5pt&\omit&height5pt&\omit&height5pt&\omit&height5pt&&\cr
&& yes && no && yes && M(6,\{e\},0) &\cr
\omit&height5pt&\omit&height5pt&\omit&height5pt&\omit&height5pt&&\cr\tablerule
\omit&height5pt&\omit&height5pt&\omit&height5pt&\omit&height5pt&&\cr
&& yes && no && no && ? &\cr
\omit&height5pt&\omit&height5pt&\omit&height5pt&\omit&height5pt&&\cr\tablerule
\omit&height5pt&\omit&height5pt&\omit&height5pt&\omit&height5pt&&\cr
&& no && yes && yes && ? &\cr
\omit&height5pt&\omit&height5pt&\omit&height5pt&\omit&height5pt&&\cr\tablerule
\omit&height5pt&\omit&height5pt&\omit&height5pt&\omit&height5pt&&\cr
&& no && yes && no && Impossible &\cr
\omit&height5pt&\omit&height5pt&\omit&height5pt&\omit&height5pt&&\cr\tablerule
\omit&height5pt&\omit&height5pt&\omit&height5pt&\omit&height5pt&&\cr
&& no && no && yes && $\widetilde{\Bbb C P^5}\times \widetilde{\Bbb C
P^5}$ &\cr
\omit&height5pt&\omit&height5pt&\omit&height5pt&\omit&height5pt&&\cr\tablerule
\omit&height5pt&\omit&height5pt&\omit&height5pt&\omit&height5pt&&\cr
&& no && no && no && $\widetilde{\Bbb C P^5}$ &\cr
\omit&height5pt&\omit&height5pt&\omit&height5pt&\omit&height5pt&&\cr
\tablerule\noalign{\smallskip}\hfil\cr}}
\vfill
\pagebreak

\p We prove the theorem via line-by-line analysis of Tables 1 and 2.
\subhead Line 1 in Tables 1 and 2
\endsubhead
For closed K\"ahler manifolds, the Hard Lefschetz property is proved in
\cite{GH}, the evenness of $b_{2i+1}$ follows from the Hodge theory
\cite{W}, the triviality of Massey products follows from the formality
of any closed K\"ahler manifold \cite{DGMS}.

One can ask if there are non-K\"ahler manifolds having the properties
from line 1. In the symplectically aspherical case the answer is
affirmative. Let $G=\Bbb R\times _{\phi}\Bbb R^2$ be the semidirect
product determined by the one-parameter subgroup
$\phi(t)=\text{diag}\,(e^{kt},e^{-kt}), t\in\Bbb R, e^k+e^{-k}\not=2$.
One can check that $G$ contains a lattice, say $\Gamma$. Then the
compact solvmanifold
$$M=G/\Gamma\times S^1$$
is symplectic and has the same minimal model as the K\"ahler manifold
$S^2\times T^2$. Hence such manifold fits into line 1. It cannot be
K\"ahler, since it admits no complex structure. The latter follows
from the Kodaira--Yau classification of compact complex surfaces (see
\cite{TO} for details).

\subhead Line 2 in Tables 1 and 2 \endsubhead
Any manifold satisfying the Hard Lefschetz property must have {\it
even}  $b_{2i+1}$. Indeed, consider the usual non-singular pairing $p:
H^{2k+1}(M)\otimes H^{2m-2k-1}(M)\to \Bbb R$ of the form
$$
p\left([\alpha],[\beta]\right)=\int_M\alpha\wedge\beta.
$$
Define a skew-symmetric bilinear form $\langle -,-\rangle:
H^{2k+1}(M)\otimes H^{2k+1}(M)\to\Bbb R$ via the formula
$$
\langle [\alpha],[\gamma]\rangle=p\left([\alpha],L^{m-2k-1}[\gamma]\right),
$$
for $[\alpha],[\gamma]\in H^{2k+1}(M)$. Since this form is
non-degenerate and skew-symmetric, its domain $H^{2k+1}(M)$ must be
even-dimensional, i.e. $b_{2k+1}$ is even.

\subhead Line 3 in Table 1 \endsubhead
We do not know any non-simply-connected (and, in particular,
symplectically aspherical) examples to fill in this line.

\subhead Line 3 in Table 2 \endsubhead
We use Theorem 2.2 with $n=6$ and $G=\{e\}$. Then, for every $b$, the
corresponding manifold $M(6,\{e\},b)$ has even odd-degree Betti
numbers and does not have the Hard Lefschetz property. Furthermore,
all the Massey products in $M$ vanish by 2.1.

\subhead Line 4 and 5 in Tables 1 and 2 \endsubhead
We do not know any examples to fill in these lines.

\subhead Line 6 in Tables 1 and 2 \endsubhead
This is impossible, see the argument concerning line 2.

\subhead Lines 7 and 8 in Table 1\endsubhead
Consider the Kodaira-Thurston manifold $K$ \cite{Th}. Recall that this
manifold is defined as a nilmanifold
$$
K=N_3/\Gamma\times S^1,
$$
where $N_3$ denotes the 3-dimensional nilpotent Lie group of
triangular unipotent matrices and $\Gamma$ denotes the lattice of such
matrices with integer entries. One can check that the
Chevalley--Eilenberg complex of the Lie algebra $\frak n_3$ is of the
form
$$(\Lambda(e_1,e_2,e_3),d), \,\, de_1=de_2=0,\, de_3=e_1e_2.$$
with $|e_i|=1$. We have already mentioned that the minimal model of
any nilmanifold $N/\Gamma$ is isomorphic to the Chevalley--Eilenberg
complex of the Lie algebra $\frak n$. In particular, one can get the
minimal model of the Kodaira--Thurston manifold in the form
$$
(\Lambda(x, e_1, e_2, e_3), d),\,\,\ dx=de_1=de_2=0,\,\, de_3=e_1e_2
$$
with degrees of all generators equal 1. One can check that the vector
space $H^1(K)$ has the basis \{$[x], [e_1], [e_2]\}$. Hence,
$b_1(K)=3$, which also shows that $K$ does not have the Hard Lefschetz
property. Furthermore, $K$ possesses a symplectic form $\omega$ with
$[\omega]=[e_1e_3+e_2x]$, and one can prove that the Massey triple
product $\langle [e_1],[e_1],[\omega]\rangle$ is non-trivial. Thus,
$K$ realizes Line 8 of Table 1.

Finally, $K \times K$ realizes Line 7 of Table 1.

\subhead Lines 7 and 8 in Table 2 \endsubhead
We use Theorem 2.5. Consider a symplectic embedding $i: K \to \Bbb CP^m,
m\geq 5$, and perform the blow-up along $i$. Then, by 2.5(i),
$\widetilde{\Bbb CP^m}$ is simply-connected. Furthermore, it realizes Line 8 of Table 2 by 2.5(ii) and 2.5(iii).

Finally, $\widetilde{\Bbb CP^m} \times \widetilde{\Bbb CP^m}$ realizes
Line 7 of Table 2.
\qed
\m
{\noindent \bf 3.2. Remark.} The result of Lupton \cite{L} shows that
the problem of constructing of a non-formal  manifold with the Hard
Lefschetz property turns our to be very delicate. In \cite{L} there is
an example of a DGA, whose cohomology has the Hard Lefschetz property,
but which is not {\it intrinsically} formal. This means that there is
also a {\it non-formal} minimal algebra with the same cohomology ring.
Sometimes, using Browder--Novikov theory, one can construct a smooth
closed manifold $M$ with such non-formal Sullivan minimal model.
However, there is no way in sight to get a symplectic structure on
$M$.

\head 4. Flexible symplectic manifolds \endhead

Let $(M^{2m},\omega)$ be a symplectic manifold. It is known that there
exists a unique non-degenerate Poisson structure $\Pi$ associated with
the symplectic structure (see, for example \cite{LM, TO}). Recall
that $\Pi$ is a skew symmetric tensor field of order 2 such that
$[\Pi,\Pi]=0$, where $[-,-]$ is the Schouten-Nijenhuis bracket.
\m
The Koszul differential $\delta: \Omega^k(M)\to \Omega^{k-1}(M)$ is
defined for Poisson, in particular symplectic, manifolds as
$$
\delta=[i(\Pi),d].
$$
Brylinski has proved in \cite{B} that the Koszul differential is a
symplectic codifferential of the exterior differential with respect to
the symplectic star operator. We choose the volume form associated to
the symplectic form, say $v_M=\omega^m/m!$. Then we define the
symplectic star operator
$$
*: \Omega^k(M)\to\Omega^{2m-k}(M)
$$
by the condition $\beta\wedge (*\alpha)=\Lambda^k(\Pi)(\beta,
\alpha)v_M$, for all $\alpha,
\beta\in \Omega^k(M)$. It turns out to be that
$$
\delta=(-1)^{k+1}(*\circ d\circ *).
$$
\m
{\noindent \bf 4.1. Definition.} A $k$-form $\alpha$ on the symplectic
manifold $M$ is called {\it symplectically harmonic}, if
$d\alpha=0=\delta\alpha$.

We denote by $\Omega_{\text{hr}}^k(M)$ the space of symplectically
harmonic $k$-forms on $M$. We set
$$
H^k_{\text{hr}}(M, \omega)=\Omega_{\text{hr}}^k(M)/(\Im
d\cap\Omega^k_{\text{hr}}(M)),\quad 
h_k(M)=h_k(M,\omega)=\dim\,H^k_{\text{hr}}(M,\omega).
$$
We say that a de Rham cohomology class is {\it symplectically
harmonic} if it contains a symplectically harmonic representative,
i.e. if it belongs to the subgroup $H^*_{\text{hr}}(M)$ of $H^*(M)$.
\m
{\noindent \bf 4.2. Definition.} We say that a closed smooth manifold
$M$ is {\it flexible}, if $M$ possesses a continuous family of
symplectic forms $\omega_t, t\in [a,b]$, such that
$h_k(M,\omega_a)\not=h_k(M,\omega_b)$ for some $k$.

So, the McDuff--Khesin Question (see the introduction) asks about existence of flexible
manifolds.
\m
In order to prove our result on the existence of flexible
6-dimensional nilmanifolds, we need some preliminaries. The following
lemma is proved in \cite{IRTU} and generalizes an observation of
Yan~\cite{Y}.
\proclaim{4.3. Lemma} For any symplectic manifold $(M^{2m},\omega)$ and $k=0,1,2$ we have
$$
H^{2m-k}_{\text{hr}}(M)=\Im \{L^{m-k}: H^k(M)\to H^{2m-k}(M)\}\subset
H^{2m-k}(M).
$$
In other words,
$$
h_{2m-k}(M,\omega)= \dim \Im \{L^{m-k}: H^k(M)\to H^{2m-k}(M)\}.\quad \qed
$$
\endproclaim
The following fact can be deduced from 4.3 using standard arguments
from linear algebra, see \cite{IRTU}.
\proclaim{4.4. Lemma} Let $\omega_1$ and $\omega_2$ be two symplectic
forms on a closed manifold $M^{2m}$. Suppose that, for $k=1$ or $k=2$, we have
$$
h_{2m-k}(M,\omega_1) \neq h_{2m-k}(M,\omega_2).
$$
Then $M$ is flexible.
\qed
\endproclaim
\proclaim{4.5. Proposition} Let $G$ be a simply connected $6$-dimensional
nilpotent Lie group such that its Lie algebra $\frak g$ has the basis
$\{X_i\}_{i=1}^6$ and the following structure relations:
$$
 [X_1,X_2] = - X_4,\quad  
    [X_1,X_4] = - X_5,\quad 
    [X_1, X_5] = [X_2,X_3] = [X_2,X_4] = - X_6
$$
$($all the other brackets $[X_i,X_j]$ are assumed to be zero$)$. Then
$G$ admits a lattice $\Gamma$, and the corresponding compact
nilmanifold $N:=G/\Gamma$ admits two symplectic forms $\omega_1$ and
$\omega_2$ such that
$$
\dim \Im L^2_{[\omega_1]}=0,\quad \dim \Im L^2_{[\omega_2]}=2.
$$
\endproclaim
\p First, $G$ has a lattice by 2.4(ii). Furthermore, by 2.4(iii), in the
Chevalley--Eilenberg complex $(\Lambda^*\frak g^*,d)$ we have
$$
\split d\alpha_1&=d\alpha_2=d\alpha_3=0,\\
d\alpha_4&=\alpha_1\alpha_2,\\
d\alpha_5&=\alpha_1\alpha_4,\\
d\alpha_6&=\alpha_1\alpha_5+\alpha_2\alpha_3+\alpha_2\alpha_4,
\endsplit
$$
where we write $\alpha_i\alpha_j$ instead of $\alpha_i\wedge \alpha_j$.
One can check that the following elements represent closed homogeneous
2-forms on $N$:
$$
\split
\omega_1&=\alpha_1\alpha_6+\alpha_2\alpha_5-\alpha_3\alpha_4,\\
\omega_2&=\alpha_1\alpha_3+\alpha_2\alpha_6-\alpha_4\alpha_5.
\endsplit
$$
Since $[\omega^3_1]\neq 0\neq [\omega_2^3]$, these homogeneous forms
are symplectic. Indeed, by 2.4(iii) the cohomology classes
$[\omega_0]$ and $[\omega_1]$ have homogeneous representatives whose
third powers are non-zero. Then the same is valid for their pull-backs
to {\it invariant} 2-forms on the Lie group $G$. But for invariant
2-forms this condition implies non-degeneracy. Since $G \to N$ is a
covering, the homogeneous forms $\omega_1$ and $\omega_2$ on $N$ are
also non-degenerate. 
\m 
Obviously, the $\Bbb R$-vector space $H^1(N)$ has the basis $\{[\alpha_1],[\alpha_2],[\alpha_3]\}$.
One can check by direct calculation that
$$
[\omega_1]^2[\alpha_i]=0, \,\,i=1,2,3
$$
and that
$$
[\omega_2]^2[\alpha_1]=-2[\alpha_1\alpha_2\alpha_4\alpha_5\alpha_6],
\,\,
[\omega_2]^2[\alpha_2]=0,\,\,[\omega_2]^2[\alpha_3]=2[\alpha_2\alpha_3\alpha_4\alpha_5\alpha_6].
$$
Finally, it is straightforward that the above cohomology classes span
2-dimensional subspace in $H^5(N)$.
\qed
\proclaim{4.6. Theorem} There exists a flexible $6$-dimensional nilmanifold.
\endproclaim
\p
Consider the nilmanifold $N$ as in 4.5. Because of 4.3 and 4.5, we conclude that
$$
h_5(N,\omega_1)=0\neq 2=h_5(N,\omega_2),
$$
and the result follows from 4.4.
\qed

{\bf Acknowledgment.} The first and the fourth authors were partially supported by the project UPV 127.310-EA147/98. This work was partially done in Oberwolfach and financed by Volkswagen-Stiftung. The second and third authors were also partially supported by Max-Planck  Institut f\"ur Mathematik, Bonn.

\enddocument
\end